\documentclass[11pt]{article}

\usepackage[in]{fullpage}
\usepackage{times}
\usepackage{graphicx}
\usepackage[rflt]{floatflt}
\usepackage{epsfig,subfigure,color,multirow}
\usepackage{amsmath,amssymb,algorithm,algorithmic,theorem,float,bbm,bm,enumerate}
\usepackage{url}
\usepackage{epstopdf}

\newcommand{\beq}{\begin{equation}}
\newcommand{\eeq}{\end{equation}}

\newcommand\R{\mathbb{R}}

\newcommand{\bd}{\boldsymbol}

\renewcommand{\a}{\mathbf{a}}
\renewcommand{\b}{\mathbf{b}}
\renewcommand{\c}{\mathbf{c}}
\newcommand{\e}{\mathbf{e}}

\renewcommand{\u}{\mathbf{u}}
\renewcommand{\v}{\mathbf{v}}
\newcommand{\w}{\mathbf{w}}
\newcommand{\x}{\mathbf{x}}
\newcommand{\y}{\mathbf{y}}
\newcommand{\z}{\mathbf{z}}

\newcommand{\cX}{{\cal X}}

\newcommand{\cZ}{{\cal Z}}

\newcommand{\bA}{\mathbf{A}}
\newcommand{\bB}{\mathbf{B}}

\newcommand{\bC}{\mathbf{C}}

\newcommand{\bZ}{\mathbf{Z}}

\newcommand{\bX}{\mathbf{X}}
\newcommand{\bY}{\mathbf{Y}}

\newcommand{\bI}{\mathbf{I}}

\newcommand{\myref}[1]{(\ref{#1})}

\DeclareMathOperator{\argmin}{argmin}

 \DeclareMathOperator{\tr}{Tr}

\newcounter{exampleI}
{\theorembodyfont{\rmfamily} \newtheorem{exm}{Example}}
\newtheorem{thm}{Theorem}
\newtheorem{prop}{Proposition}
\newtheorem{lem}{Lemma}
\newtheorem{asm}{Assumption}
\newtheorem{cor}{Corollary}

\newtheorem{remark}{Remark}
\newcommand{\proof}{\noindent{\itshape Proof:}\hspace*{1em}}

\newcommand{\qed}{\nolinebreak[1]~~~\hspace*{\fill} \rule{5pt}{5pt}\vspace*{\parskip}\vspace*{1ex}}

\newcommand {\commentout}[1] {}

\title{Bregman Alternating Direction Method of Multipliers}
\author{Huahua Wang \\
Dept of Computer Science \& Engg\\
University of Minnesota, Twin Cities\\
huwang@cs.umn.edu
\and
Arindam Banerjee\\
Dept of Computer Science \& Engg\\
University of Minnesota, Twin Cities\\
banerjee@cs.umn.edu
}

\date{}

\begin{document}

\maketitle

\begin{abstract}
The mirror descent algorithm (MDA) generalizes gradient descent by using a Bregman divergence to replace squared Euclidean distance. In this paper,
we similarly generalize the
alternating direction method of multipliers (ADMM) to Bregman ADMM (BADMM), which allows the choice of different Bregman divergences to exploit the structure of problems. BADMM provides a unified framework for ADMM and its variants, including generalized ADMM, inexact ADMM and Bethe ADMM.
 We establish the global convergence and the $O(1/T)$ iteration complexity for BADMM. In some cases,  BADMM can be faster than ADMM by a factor of $O(n/\log(n))$ where $n$ is the dimension of the problem. Experimental results are illustrated on the mass transportation problem, which can be solved in parallel by BADMM. BADMM is faster than ADMM and highly optimized commercial software Gurobi, particularly when implemented on GPU.
\end{abstract}

\section{Introduction}\label{sec:intro}
In recent years, the alternating direction method of multipliers (ADM or ADMM)~\cite{boyd10} has been successfully applied
in a broad spectrum of applications, ranging from image processing~\cite{Figueiredo10,osher10:admm} to applied statistics and machine learning~\cite{yang09,wang12:oadm,wang13:ADMMMAP}. 
For further understanding of ADMM, we refer the readers to the comprehensive review by~\cite{boyd10} and references therein.
In particular, ADMM considers the problem of minimizing composite objective functions subject to an equality constraint:
\begin{align} \label{eq:adm_pm}
\min_{\x\in \mathcal{X},\z\in\mathcal{Z}} f(\x) + g(\z) \quad \text{s.t.} \quad \bA\x + \bB \z = \c~,
\end{align}
where $f$ and $g$ are convex functions, $\bA \in \R^{m\times n_1}, \bB \in \R^{m\times n_2}, \c \in \R^{m\times 1}$, $\x \in \mathcal{X} \in \R^{n_1\times 1}, \z \in \mathcal{Z} \in \R^{n_2\times 1}$, and $\mathcal{X}$ and $\mathcal{Z}$ are convex sets. $f$ and $g$ can be non-smooth functions, including indicator functions of convex sets. Many machine learning problems can be cast into the framework of minimizing a composite objective~\cite{nest07:composite,Duchi10_comid}, where $f$ is a loss function such as hinge or logistic loss, and $g$ is a regularizer, e.g., $\ell_1$ norm, $\ell_2$ norm, nuclear norm or total variation. The two functions usually have different structures and constraints because they have different tasks in data mining. Therefore, it is useful and sometimes necessary to split and solve them separately, which is exactly the forte of ADMM.

In each iteration, ADMM updates splitting variables separately and alternatively by solving the augmented Lagrangian of~\myref{eq:adm_pm}, which is defined as follows:
\begin{align}\label{eq:lag1}
L_{\rho}(\x,\z,\y) & = f(\x) + g(\z) + \langle \y , \bA\x + \bB\z - \c \rangle + \frac{\rho}{2} \| \bA\x + \bB\z - \c \|_2^2,
\end{align}
where $\y \in \R^m$ is dual variable, $\rho > 0 $ is penalty parameter, and the quadratic penalty term is to penalize the violation of the equality constraint. ADMM consists of the following three updates:
\begin{align}
&\x_{t+1} = \argmin_\x ~f(\x) + \langle \y_t , \bA\x + \bB\z_t - \c \rangle + \frac{\rho}{2} \| \bA\x + \bB\z_t - \c \|_2^2 ~,  \label{eq:adm_x}\\
&\z_{t+1}  = \argmin_\z~g(\z) + \langle \y_t , \bA\x_{t+1} + \bB\z - \c \rangle + \frac{\rho}{2} \| \bA\x_{t+1} + \bB\z - \c \|_2^2~,  \label{eq:adm_z} \\
&\y_{t+1}  = \y_t + \rho(\bA\x_{t+1} + \bB\z_{t+1} - \c)~. \label{eq:adm_y}
\end{align}
Since the computational complexity of $\y$ update~\myref{eq:adm_y} is trivial, the computational complexity of ADMM lies in the $\x$ and $\z$ updates~\myref{eq:adm_x}-\myref{eq:adm_z} which amount to solving proximal minimization problems using the quadratic penalty term. 
Inexact ADMM~\cite{yang09,boyd10} and generalized ADMM~\cite{deng12:admm} have also been proposed to solve the updates inexactly by linearizing the functions and adding additional quadratic terms. Recently, online ADMM~\cite{wang12:oadm} and Bethe-ADMM~\cite{wang13:ADMMMAP} add an additional Bregman divergence on the $\x$ update by keeping or linearizing the quadratic penalty term $\| \bA\x + \bB\z - \c \|_2^2$. As far as we know, all existing ADMMs use quadratic penalty terms.

A large amount of literature shows that replacing the quadratic term by Bregman divergence in gradient-type methods could greatly boost their performance in solving the constrained optimization problem. 
First, the use of Bregman divergence could effectively exploit the structure of problems~\cite{chen93:proxBreg,Beck03,Duchi10_comid} , e.g.,  in computerized tomography~\cite{ben01:mda}, clustering problems and exponential family distributions~\cite{bane05:bregman}. 
Second, in some cases, the gradient descent method with Kullback-Leibler (KL) divergence can outperform the method with the quadratic term by a factor of $O(\sqrt{n\ln n})$ where $n$ is the dimensionality of the problem~\cite{Beck03,ben01:mda}. 
Mirror descent algorithm (MDA) and composite objective mirror descent (COMID)~\cite{Duchi10_comid} use Bregman divergence to replace the quadratic term in gradient descent or proximal gradient~\cite{comb09:prox}. 
Proximal point method with D-functions (PMD)~\cite{chen93:proxBreg,ceze98} and Bregman proximal minimization (BPM) ~\cite{kiwiel95:gbregman} generalize proximal point method by using Bregman divegence to replace the quadratic term. 

On the side of ADMM, it is still unknown whether the quadratic penalty term in ADMM can be replaced by Bregman divergence, although the convergence of ADMM is well understood. 
The proof of global convergence of ADMM can be found in \cite{Gabay83,boyd10}. Recently, it has been shown that ADMM converges at a rate of $O(1/T)$~\cite{wang12:oadm,he12:vi}, where $T$ is the number of iterations. For strongly convex functions, the dual objective of an accelerated version of ADMM can converge at a rate of $O(1/T^2)$~\cite{Goldstein12:fadmm}. Under suitable assumptions like strongly convex functions or a sufficiently small step size for the dual variable update, ADMM can achieve a linear convergence rate~\cite{deng12:admm,luo12:admm}.
However, as pointed out by~\cite{boyd10}, \lq\lq There is currently no proof of convergence known for ADMM with nonquadratic penalty terms.\rq\rq


In this paper, we propose Bregman ADMM (BADMM) which uses Bregman divergences to replace the quadratic penalty term in ADMM, answering the question raised in~\cite{boyd10}.
More specifically, the quadratic penalty term in the $\x$ and $\z$ updates~\myref{eq:adm_x}-\myref{eq:adm_z} will be replaced by a Bregman divergence in BADMM. We also introduce a generalized version of BADMM where two additional Bregman divergences are added to the $\x$ and $\z$ updates. The generalized BADMM (BADMM for short) provides a unified framework for solving~\myref{eq:adm_pm}, which allows one to choose suitable Bregman divergence so that the $\x$ and $\z$ updates can be solved efficiently. BADMM includes ADMM and its variants as special cases. In particular, BADMM replaces all quadratic terms in generalized ADMM~\cite{deng12:admm} with Bregman divergences. By choosing a proper Bregman divergence, we also show that inexact ADMM~\cite{yang09} and Bethe ADMM~\cite{wang13:ADMMMAP} can be considered as special cases of BADMM. BADMM generalizes ADMM similar to how MDA generalizes gradient descent and how PMD generalizes proximal methods. In BADMM, the $\x$ and $\z$ updates can take the form of MDA or PMD.
We establish the global convergence and the $O(1/T)$ iteration complexity for BADMM. In some cases, we show that BADMM can outperform ADMM by a factor $O(n/\ln n)$. 
We evaluate the performance of BADMM in solving the linear program problem of mass transportation~\cite{cock41:mt}.  By exploiting the structure of the problem, BADMM leads to massive parallelism and can easily run on GPU. BADMM can even be orders of magnitude faster than highly optimized commercial software Gurobi. While Gurobi breaks down in solving a linear program of hundreds of millions of parameters in a server, BADMM takes hundreds of seconds running in a single GPU.

The rest of the paper is organized as follows. In Section 2, we propose Bregman ADMM and discuss several special cases of BADMM. In Section 3, we establish the convergence of BADMM. In Section 4, we consider illustrative applications of BADMM, and conclude in Section 5.

\section{Bregman Alternating Direction Method of Multipliers}\label{sec:badm}


Let $\phi:\Omega\rightarrow\R$ be a continuously differentiable and strictly convex function on the relative interior of a convex set $\Omega$. Denote $\nabla\phi(\y)$ as the gradient of $\phi$ at $\y$. We define Bregman divergence\footnote{The definition of Bregman divergence has been generalized to nondifferentiable functions~\cite{kiwiel95:gbregman,teda12:gbregman}.} $B_{\phi}:\Omega\times \text{ri}(\Omega) \rightarrow \R_+$ induced by $\phi$ as
\begin{align*}
B_{\phi}(\x,\y) = \phi(\x) - \phi(\y) - \langle \nabla\phi(\y), \x - \y \rangle~.
\end{align*}
Since $\phi$ is convex, $B_{\phi}(\x,\y) \geq 0$ where the equality holds if and only if $\x = \y$. More details about Bregman divergence can be found in~\cite{chen93:proxBreg,bane05:bregman}. 
Two of the most commonly used examples are squared Euclidean distance $B_{\phi}(\x,\y) = \frac{1}{2}\|\x - \y\|_2^2$ and KL divergence $B_{\phi}(\x,\y) = \sum_{i=1}^n x_i \log \frac{x_i}{y_i}$. 

Assuming $B_{\phi} (\c - \bA\x, \bB\z)$ is well defined, we replace the quadratic penalty term in the augmented Lagrangian~\myref{eq:lag1} by a Bregman divergence as follows:
\begin{align}\label{eq:blag}
L_{\rho}^{\phi}(\x,\z,\y) &= f(\x) + g(\z) + \langle \y , \bA\x + \bB\z - \c \rangle + \rho B_{\phi} (\c - \bA\x, \bB\z).
\end{align}
Unfortunately, we can not derive Bregman ADMM (BADMM) updates by simply solving $L_{\rho}^{\phi}(\x,\z,\y)$ alternatingly as ADMM does because Bregman divergences are not necessarily convex in the second argument.
More specifically, given $(\z_t,\y_t)$, $\x_{t+1}$ can be obtained by solving $\min_{\x}~L_{\rho}^{\phi}(\x,\z_t,\y_t)$, where the quadratic penalty term $\frac{1}{2}\| \bA\x + \bB\z_t - \c \|_2^2$ for ADMM in~\myref{eq:adm_x} is replaced with $B_{\phi}(\c - \bA\x, \bB\z_t)$ in the $\x$ update of BADMM. However, given $(\x_{t+1},\y_t)$, we cannot obtain $\z_{t+1}$ by solving $\min_{\z}~L_{\rho}^{\phi}(\x_{t+1},\z,\y_t)$, since the term $B_{\phi} (\c - \bA\x_{t+1}, \bB\z)$ need not be convex in $\z$. The observation motivates a closer look at the role of the quadratic term in ADMM.

In standard ADMM, the quadratic augmentation term added to the Lagrangian is just a penalty term to ensure the new updates do not violate the constraint significantly. Staying with these goals, we propose the $\z$ update augmentation term of BADMM to be: $B_{\phi}(\bB\z, \c - \bA\x_{t+1})$, instead of the quadratic penalty term $\frac{1}{2}\| \bA\x_{t+1} + \bB\z - \c \|_2^2$ in~\myref{eq:adm_x}.
Then, we get the following updates for BADMM:
\begin{align}
\x_{t+1} = & \underset{\x\in\mathcal{X}}{\argmin}~f(\x) + \langle \y_t, \bA\x + \bB\z_t - \c \rangle  + \rho B_{\phi}(\c - \bA\x, \bB\z_t) ~,  \label{eq:ebadmm_x}\\
\z_{t+1}  = & \underset{\z\in\mathcal{Z}}{\argmin}~g(\z) + \langle \y_t, \bA\x_{t+1} + \bB\z - \c \rangle + \rho B_{\phi}(\bB\z, \c - \bA\x_{t+1})~,  \label{eq:ebadmm_z} \\
\y_{t+1}  = & \y_t + \rho(\bA\x_{t+1} + \bB\z_{t+1} - \c)~. \label{eq:ebadmm_y}
\end{align}
Compared to ADMM~\myref{eq:adm_x}-\myref{eq:adm_y}, BADMM simply uses a Bregman divergence to replace the quadratic penalty term in the $\x$ and $\z$ updates. It is worth noting that the same Bregman divergence $B_{\phi}$ is used in the $\x$ and $\z$ updates. 

We consider a special case when $\bA = -\bI, \bB = \bI, \c = \mathbf{0}$.~\myref{eq:ebadmm_x} is reduced to
\begin{align}\label{eq:snr12}
\x_{t+1} = \underset{\x\in\mathcal{X}}{\argmin}~f(\x) + \langle \y_t, - \x + \z_t \rangle + \rho B_{\phi}(\x, \z_t)~. 
\end{align}
If $\phi$ is a quadratic function, the constrained problem~\myref{eq:snr12} requires the projection onto the constraint set $\cX$. However, in some cases, if choosing a proper Bregman divergence,~\myref{eq:snr12} can be solved efficiently or has a closed-form solution. For example, if $f$ is a linear function and $\mathcal{X}$ is the unit simplex, $B_{\phi}$ should be KL divergence, leading to the exponentiated gradient~\cite{Beck03,ben01:mda,Nemi83:complexity}. 
Interestingly, if the $\z$ update is also the exponentiated gradient, we have alternating exponentiated gradients. In Section 4, we will show the mass transportation problem can be cast into this scenario. 


 While the updates~\myref{eq:ebadmm_x}-\myref{eq:ebadmm_z} use the same Bregman divergences, efficiently solving the $\x$ and $\z$ updates may not be feasible, especially when the structure of the original functions $f,g$, the function $\phi$ used for augmentation, and the constraint sets $\cX,\cZ$ are rather different. For example, if $f(\x)$ is a logistic function in ~\myref{eq:snr12}, it will not have a closed-form solution even $B_{\phi}$ is the KL divergence and $\cX$ is the unit simplex. To address such concerns, we propose a generalized version of BADMM in Section 2.1. 

\subsection{Generalized BADMM}
\label{ssec:gbamm}

To allow the use of different Bregman divergences in the $\x$ and $\z$ updates~\myref{eq:ebadmm_x}-\myref{eq:ebadmm_y} of BADMM, the generalized BADMM simply introduces an additional Bregman divergence for each update. The generalized BADMM has the following updates:
\begin{align}
\x_{t+1} = & ~\underset{\x\in\mathcal{X}}{\argmin}~f(\x) + \langle \y_t, \bA\x + \bB\z_t - \c \rangle + \rho B_{\phi}(\c - \bA\x, \bB\z_t) + \rho_{\x} B_{\varphi_{\x}}(\x,\x_t)~,  \label{eq:badmm_x}\\
\z_{t+1}  = & ~\underset{\z\in\mathcal{Z}}{\argmin}~g(\z) + \langle \y_t, \bA\x_{t+1} + \bB\z - \c \rangle + \rho B_{\phi}(\bB\z, \c - \bA\x_{t+1}) + \rho_{\z} B_{\varphi_{\z}}(\z,\z_t)~,  \label{eq:badmm_z} \\
\y_{t+1}  = & ~\y_t + \tau(\bA\x_{t+1} + \bB\z_{t+1} - \c)~. \label{eq:badmm_y}
\end{align}
where $\rho > 0, \tau > 0, \rho_\x \geq 0,  \rho_\z \geq 0$. Note that we allow the use of a different step size $\tau$ in the dual variable update~\cite{deng12:admm,luo12:admm}. There are three Bregman divergences in the generalized BADMM. While the Bregman divergence $B_{\phi}$ is shared by the $\x$ and $\z$ updates, the $\x$ update has its own Bregman divergence $B_{\varphi_{\x}}$ and the $\z$ update has its own Bregman divergence $B_{\varphi_{\z}}$. The two additional Bregman divergences in generalized BADMM are variable specific, and can be chosen to make sure that the $\x_{t+1},\z_{t+1}$ updates are efficient. If all three Bregman divergences are quadratic functions, the generalized BADMM reduces to the generalized ADMM~\cite{deng12:admm}. We prove convergence of generalized BADMM in Section 3, which yields the convergence of BADMM with $\rho_x = \rho_z = 0$.



In the following, we illustrate how to choose a proper Bregman divergence $B_{\varphi_{\x}}$ so that the $\x$ update can be solved efficiently, e.g., a closed-form solution, noting that the same arguments apply to the $\z$-updates. Consider the first three terms in~\myref{eq:badmm_x} as $s(\x) + h(\x)$, where $s(\x)$ denotes an easy term and $h(\x)$ is the problematic term which needs to be linearized for an efficient $\x$-update. We illustrate the idea with several examples later in the section. Now, we have
 \begin{align}
\x_{t+1} = \min_{\x \in \cX}~s(\x) + h(\x) + \rho_{\x}B_{\varphi_\x}(\x,\x_t)~.
 \end{align}
 where efficient updates are difficult due to the mismatch in structure between $h$ and $\cX$. The goal is to `linearize' the function $h$ by using the fact that the Bregman divergence $B_h(\x,\x_t)$ captures all the higher-order (beyond linear) terms in $h(\x)$ so that:
\beq
h(\x) - B_h(\x,\x_t) = h(\x_t) + \langle \x - \x_t, \nabla h(\x_t) \rangle~
\eeq
is a linear function of $\x$.
Let $\psi$ be another convex function such that one can efficiently solve $\min_{x \in \cX}~s(\x) + \psi(x) + \langle \x, \b \rangle$ for any constant $\b$. Assuming $\varphi_\x(\x) = \psi(\x) - \frac{1}{\eta} h(\x)$ is convex, we construct a Bregman divergence based proximal term to the original problem so that:
\begin{align} \label{eq:x_efficient}
\!\!\!\!\underset{\x \in \cX}{\argmin}~ s(\x) + h(\x)  + \rho_\x B_{\varphi_\x}(\x,\x_t) =\underset{\x \in \cX}{\argmin}~ s(\x) + \psi(\x) + \langle \x, \frac{1}{\rho_\x} \nabla h(\x_t) - \nabla \psi(\x_t) \rangle~,
\end{align}
where the latter problem can be solved efficiently, by our assumption. To ensure $\varphi_\x$ is convex, we need the following condition:
\begin{prop}
If $h$ is smooth and has Lipschitz continuous gradients with constant $\nu$ under a $p$-norm, then $\varphi_\x$ is $\nu/\rho_\x$-strongly convex w.r.t. the $p$-norm.
\end{prop}


This condition has been widely used in gradient-type methods, including MDA and COMID. Note that the convergence analysis of generalized ADMM in Section 4 holds for any additional Bregman divergence based proximal terms, and does not rely on such specific choices. Using the above idea, one can `linearize' different parts of the $\x$ update to yield an efficient update.

We consider three special cases, respectively focusing on linearizing the function $f(\x)$, linearizing the Bregman divergence based augmentation term $B_{\phi}(\c - \bA \x,\bB \z_t)$, and linearizing both terms, along with examples for each case. 

\textbf{Case 1: Linearization of smooth function $f$:}
Let $h(\x) = f(\x)$ in~\myref{eq:x_efficient}, we have
\begin{align}
\x_{t+1} &= \underset{\x\in\mathcal{X}}{\argmin}~ \langle \nabla f(\x_t), \x - \x_t\rangle + \langle \y_t, \bA\x \rangle + \rho B_{\phi}(\c - \bA\x, \bB\z_t) + \rho_\x B_{\psi_\x}(\x, \x_t)~.
\end{align}
where $\nabla f(\x_t)$ is the gradient of $f(\x)$ at $\x_t$. 
\begin{exm}\label{eg:logistic}
Consider the following ADMM form for sparse logistic regression problem~\cite{hatf09,boyd10}:
\begin{align}\label{eq:slog}
\min_\x h(\x) + \lambda \|\z\|_1~,~\text{s.t.}~ \x = \z~,
\end{align}
where $h(\x)$ is the logistic function. If we use ADMM to solve~\myref{eq:slog}, the $\x$ update is as follows~\cite{boyd10}:
\begin{align}\label{eq:logistic_x}
\x_{t+1} = \underset{\x}{\argmin}~ h(\x) + \langle \y_t, \x - \z_t\rangle + \frac{\rho}{2} \|\x - \z_t\|_2^2~,
\end{align}
which is a ridge-regularized logistic regression problem and one needs an iterative algorithm like L-BFGS to solve it. Instead, if we linearize $h(\x)$ at $\x_t$ and set $B_{\psi}$ to be a quadratic function, then
\begin{align}
\x_{t+1} = \underset{\x}{\argmin} &~\langle \nabla  ~h(\x_t), \x - \x_t \rangle + \langle \y_t, \x - \z_t \rangle + \frac{\rho}{2} \|\x - \z_t\|_2^2 + \frac{\rho_\x}{2}\|\x - \x_t\|_2^2~,
\end{align}
the $\x$ update has a simple closed-form solution.
\end{exm}

\textbf{Case 2: Linearization of the quadratic penalty term:} In ADMM, $B_{\phi}(\c - \bA\x, \bB\z_t) = \frac{1}{2} \| \bA\x + \bB\z_t - \c\|_2^2$. Let $h(\x) = \frac{1}{2} \| \bA\x + \bB\z_t - \c\|_2^2$. Then $\nabla h(\x_t) = \bA^T(\bA\x_t + \bB\z_t - \c) $, we have
\begin{align}\label{eq:snr2case2}
\x_{t+1}  =  \underset{\x\in\mathcal{X}}{\argmin} &~ f(\x) + \langle \y_t + \rho (\bA\x_t + \bB\z_t - \c), \bA\x\rangle + \rho_\x B_{\psi}(\x, \x_t)~.
\end{align}
The case mainly solves the problem due to the $\bA\x$ term which makes $\x$ updates nonseparable, whereas the linearized version can be solved with separable (parallel) updates.
Several problems have been benefited from the linearization of quadratic term~\cite{deng12:admm}, e.g., when $f$ is $\ell_1$ loss function~\cite{hatf09}, and projection onto the unit simplex or $\ell_1$ ball~\cite{Duchi08}.

\textbf{Case 3: Mirror Descent:}  In some settings, we want to linearize both the function $f$ and the quadratic augmentation term $B_{\phi}(\c - \bA\x, \bB\z_t) = \frac{1}{2} \| \bA\x + \bB\z_t - \c\|_2^2$. Let $h(\x) = f(\x) + \langle \y^t, \bA\x \rangle + \frac{\rho}{2} \| \bA\x + \bB\z_t - \c\|_2^2 $, we have
\begin{align}\label{eq:snr2case3}
\x_{t+1} &= \underset{\x\in\mathcal{X}}{\argmin} \langle \nabla h(\x_t), \x \rangle + \rho_\x B_{\psi}(\x, \x_t)~.
\end{align}
Note that \myref{eq:snr2case3} is a MDA-type update. Further, one can do a similar exercise with a general Bregman divergence based augmentation term $B_{\phi}(\c - \bA \x, \bB \z_t)$, although there has to be a good motivation for going to this route. 
\begin{exm} \label{em:betheadmm}
{\bf [Bethe-ADMM~\cite{wang13:ADMMMAP}]}
  Given an undirected graph $G = (V,E)$, where $V$ is the vertex set and $E$ is the edge set.  Assume a random discrete variable $X_i$ associated with node $i \in V$ can take $K$ values.  In a pairwise MRF, the joint distribution of a set of discrete random variables $X = \{ X_1 , \cdots, X_n \}$ ($n$ is the number of nodes in the graph) is defined in terms of nodes and cliques~\cite{wain08:graphical}. 
Consider solving the following graph-structured problem :
\begin{align}\label{eq:graphbethe}
\min~l(\bd \mu)~~\text{s.t.}~~\bd\mu \in \mathbb{L}(G)~,
\end{align}
where $l(\bd \mu)$ is a decomposable function of $\bd \mu$  and $\mathbb{L}(G)$ is the so-called local polytope~\cite{wain08:graphical} determined by the marginalization and normalization (MN) constraints for each node and edge in the graph $G$:
\begin{align}\label{eq:mncst}
\mathbb{L}(G) =  \{ \bd\mu \geq 0~,\sum_{x_i}\mu_i(x_i) = 1~,\sum_{x_j}\mu_{ij}(x_i,x_j) = \mu_i(x_i)  \}~,
\end{align}
where $\mu_i, \mu_{ij}$ are pseudo-marginal distributions of node $i$ and edge $ij$ respectively. 
In particular, ~\myref{eq:graphbethe} serves as a LP relaxation of MAP inference probem in a pairwise MRF if $ l(\bd \mu)$ is defined as follows:
\begin{align}\label{eq:maplp_badmm}
l(\bd \mu)  =  \sum_i \sum_{x_i} \theta_i(x_i) \mu_i(x_i) + \sum_{ij \in E} \sum_{x_{ij}} \theta_{ij}(x_i,x_j) \mu_{ij}(x_i,x_j),
\end{align}
where $\theta_i, \theta_{ij}$ are the potential functions of node $i$ and edge $ij$ respectively.

The complexity of polytope $\mathbb{L}(G)$ makes~\myref{eq:graphbethe} difficult to solve. One possible way is to decompose the graph into trees such that
\begin{align}
\min~\sum_{\tau}c_{\tau}l_{\tau}(\bd \mu_{\tau})~~ \text{s.t.}~~\bd\mu_{\tau} \in \mathbb{T}_\tau, \bd \mu_{\tau} = \mathbf{m}_{\tau}~,
\end{align}
where $\mathbb{T}_\tau$ denotes the MN constraints~\myref{eq:mncst} in the tree $\tau$.  $\bd \mu_\tau$ is a vector of pseudo-marginals of nodes and edges in the tree $\tau$. $\mathbf{m}$ is a global variable  which contains all trees and $\mathbf{m}_\tau$ corresponds to the tree $\tau$ in the global variable. $c_{\tau}$ is the weight for sharing variables.  The augmented Lagrangian is 
\begin{align}
L_{\rho}(\bd\mu_{\tau}, \mathbf{m}, \bd\lambda_{\tau}) = &\sum_{\tau} c_{\tau} l_{\tau} (\bd \mu_{\tau}) + \langle \bd\lambda_{\tau}, \bd\mu_{\tau} - \mathbf{m}_{\tau}\rangle + \frac{\rho}{2} \| \bd\mu_{\tau} -\mathbf{m}_{\tau} \|_2^2~.
\end{align}
which leads to the following update for $\bd\mu^{t+1}_{\tau}$ in ADMM:
\begin{align}
&\bd\mu^{t+1}_{\tau}=\underset{\bd\mu_{\tau} \in \mathbb{T}_\tau}{\argmin}~c_{\tau}l_{\tau} (\bd \mu_{\tau}) + \langle \bd\lambda_{\tau}^t, \bd\mu_{\tau} \rangle + \frac{\rho}{2} \| \bd\mu_{\tau} -\mathbf{m}_{\tau}^t \|_2^2~ \label{eq:betheadmm_x} 
\end{align}
~\myref{eq:betheadmm_x}  is difficult to solve due to the MN constraints in the tree. Let $h(\bd\mu_{\tau})$ be the objective of~\myref{eq:betheadmm_x}. If linearizing $h(\bd\mu_{\tau})$  and adding a Bregman divergence in~\myref{eq:betheadmm_x},  we have:
\begin{align}
&\bd \mu_{\tau}^{t+1} = \underset{\bd\mu_{\tau} \in \mathbb{T}_\tau}{\argmin}~\langle \nabla h(\bd\mu_{\tau}^t), \bd \mu_{\tau} \rangle + \rho_\x B_{\psi}(\bd\mu_{\tau},\bd\mu_{\tau}^t) \nonumber \\
&= \underset{\bd\mu_{\tau} \in \mathbb{T}_\tau}{\argmin}~\langle \nabla h(\bd\mu_{\tau}^t) - \rho_\x \nabla \psi(\bd \mu_{\tau}^t), \bd \mu_{\tau} \rangle + \rho_\x\psi(\bd\mu_{\tau})~, \nonumber
\end{align}
If $\psi(\bd\mu_{\tau})$ is the negative Bethe entropy of $\bd\mu_{\tau}$, the update of $\bd \mu_{\tau}^{t+1} $ becomes the Bethe entropy problem~\cite{wain08:graphical} and can be solved exactly by the sum-product algorithm in a linear time in the tree.
\end{exm}

\section{Convergence Analysis of BADMM}
\label{sec:cvg}
We need the following assumption in establishing the convergence of BADMM:
\begin{asm}\label{asm:badmm}
\quad

(a) $f:\R^{n_1}\rightarrow\R\cup\{+\infty\}$ and $g:\R^{n_2}\rightarrow\R\cup\{+\infty\}$ are closed, proper and convex. 

(b) An optimal solution exists.

(c) The Bregman divergence $B_{\phi}$ is defined on an $\alpha$-strongly convex function $\phi$ with respect to a $p$-norm $\|\cdot\|_p^2$, i.e., $B_{\phi}(\u,\v) \geq \frac{\alpha}{2} \|\u - \v\|_p^2$, where $\alpha > 0$.
\end{asm}

We start wth the Lagrangian, which is defined as follows:
\begin{equation}\label{eq:lag}
L(\x,\y,\z) = f(\x) + g(\z) + \langle \y, \bA\x + \bB\z - \c \rangle.
\end{equation}
Assume that $\{ \x^*, \z^*,\y^*\}$ satisfies the KKT conditions of~\myref{eq:lag}, i.e.,
\begin{align}
-\bA^T\y^* & \in \partial f(\x^*)~, \label{eq:kktx}\\
-\bB^T\y^* & \in \partial g(\z^*)~,\label{eq:kktz}\\
\bA\x^* + \bB\z^* - \c & = 0~. \label{eq:kkty}
\end{align}
$\{ \x^*, \z^*,\y^*\}$ is an optimal solution. The optimality conditions of~\myref{eq:badmm_x} and~\myref{eq:badmm_z} are
\begin{align}
- \bA^T\{ \y_t + \rho ( - \nabla\phi  (\c - \bA\x_{t+1}) + \nabla\phi  (\bB\z_t)\} - \rho_{\x}(\nabla\varphi_{\x}(\x_{t+1}) - \nabla\varphi_{\x}(\x_{t})) &\in \partial f(\x_{t+1})~,\label{eq:fgd}  \\
- \bB^T\{ \y_t + \rho ( \nabla\phi  (\bB\z_{t+1}) - \nabla\phi  (\c - \bA\x_{t+1})\} - \rho_{\z}(\nabla\varphi_{\z}(\z_{t+1}) - \nabla\varphi_{\z}(\z_{t})) &\in \partial g(\z_{t+1})~.\label{eq:ggd}
\end{align}
If $\bA\x_{t+1} + \bB\z_{t+1} = \c$, then $\y_{t+1} = \y_t$. Therefore, \myref{eq:kktx} is satisfied if $\bA\x_{t+1} + \bB\z_{t} = \c~, \x_{t+1} = \x_t$ in~\myref{eq:fgd}. Similarly, \myref{eq:kktz} is satisfied if $\z_{t+1} = \z_t$ in \myref{eq:ggd}. Overall, the KKT conditions~\myref{eq:kktx}-\myref{eq:kkty} are satisfied if the following optimality conditions are satisfied:
\begin{subequations}\label{eq:optcd}
\begin{align}
& B_{\varphi_{\x}}(\x_{t+1},\x_{t}) = 0~,~B_{\varphi_{\z}}(\z_{t+1},\z_{t}) = 0~, \\
&\bA\x_{t+1} + \bB\z_{t} - \c = 0~,~ \bA\x_{t+1} + \bB\z_{t+1} - \c = 0~.\label{eq:optcd0}
\end{align}
\end{subequations}
For the exact BADMM, $\rho_\x = \rho_\z = 0$ in~\myref{eq:badmm_x} and~\myref{eq:badmm_z}, the optimality conditions are~\myref{eq:optcd0},
which is equivalent to the optimality conditions used in the proof of ADMM in~\cite{boyd10}, i.e., 
\begin{align}
\bB\z_{t+1} - \bB\z_{t} = 0~,~ \bA\x_{t+1} + \bB\z_{t+1} - \c = 0~.
\end{align}

Define the residuals of optimality conditions~\myref{eq:optcd} at $(t+1)$ as:
\begin{align}\label{eq:residual}
R(t+1) = \frac{\rho_{\x}}{\rho}B_{\varphi_{\x}}(\x_{t+1},\x_t) + \frac{\rho_{\z}}{\rho}B_{\varphi_{\z}}(\z_{t+1},\z_t) + B_{\phi}(\c -\bA\x_{t+1}, \bB\z_t) +\gamma\| \bA\x_{t+1}  + \bB\z_{t+1} - \c \|_2^2~,
\end{align}
where $\gamma > 0$. If $R(t+1) = 0$, the optimality conditions~\myref{eq:optcd} and~\myref{eq:optcd0} are satisfied.
It is sufficient to show the convergence of BADMM by showing $R(t+1)$ converges to zero. We need the following lemma.
\begin{lem}\label{lem:fgobj}
Let the sequence $\{\x_t, \z_t,\y_t\}$ be generated by Bregman ADMM~\myref{eq:badmm_x}-\myref{eq:badmm_y}. For any $\x^*,\z^*$ satisfying $\bA\x^* + \bB\z^* = \c$, we have
\begin{align}\label{eq:fg}
&f(\x_{t+1})+ g(\z_{t+1}) - (f(\x^*)+g(\z^*)) \nonumber \\
&\leq -\langle \y_t, \bA\x_{t+1} + \bB\z_{t+1} - \c \rangle - \rho ( B_{\phi}(\c -\bA\x_{t+1}, \bB\z_t) +  B_{\phi}(\bB\z_{t+1}, \c -\bA\x_{t+1}))\nonumber \\
&+ \rho (B_{\phi}(\bB\z^*, \bB\z_t) - B_{\phi}(\bB\z^*, \bB\z_{t+1}))+\rho_{\x}(B_{\varphi_{\x}}(\x^*, \x_t) - B_{\varphi_{\x}}(\x^*, \x_{t+1}) - B_{\varphi_{\x}}(\x_{t+1},\x_t))\nonumber \\
& + \rho_{\z}(B_{\varphi_{\z}}(\z^*, \z_t) - B_{\varphi_{\z}}(\z^*, \z_{t+1}) - B_{\varphi_{\z}}(\z_{t+1},\z_t))
~.
\end{align}
\end{lem}
\proof
Using the convexity of $f$ and its subgradient given in~\myref{eq:fgd}, we have
\begin{align}
&f(\x_{t+1})  - f(\x) \nonumber \\
&\leq \langle - \bA^T\{ \y_t + \rho ( - \nabla\phi  (\c - \bA\x_{t+1}) + \nabla\phi  (\bB\z_t)\}- \rho_{\x}(\nabla\varphi_{\x}(\x_{t+1}) - \nabla\varphi_{\x}(\x_{t})), \x_{t+1} - \x \rangle \nonumber \\
& = -\langle \y_t, \bA(\x_{t+1} - \x) \rangle + \rho  \langle \nabla\phi  (\c - \bA\x_{t+1}) - \nabla\phi  (\bB\z_t), \bA(\x_{t+1}-\x) \rangle \nonumber \\
&- \rho_{\x}\langle \nabla\varphi_{\x}(\x_{t+1}) - \nabla\varphi_{\x}(\x_{t}), \x_{t+1} - \x \rangle~.
\end{align}
Setting $\x = \x^*$ and using $\bA\x^* + \bB\z^* = \c$, we have
\begin{align}\label{eq:fconv}
&f(\x_{t+1})  - f(\x^*) \nonumber \\
&\leq -\langle \y_t, \bA\x_{t+1} + \bB\z^* -\c \rangle + \rho  \langle \nabla\phi(\c - \bA\x_{t+1}) - \nabla\phi  (\bB\z_t), \bB\z^* - (\c - \bA\x_{t+1}) \rangle \nonumber \\
&- \rho_{\x}\langle \nabla\varphi_{\x}(\x_{t+1}) - \nabla\varphi_{\x}(\x_{t}), \x_{t+1} - \x \rangle \nonumber \\
& = -\langle \y_t, \bA\x_{t+1} + \bB\z^* -\c \rangle + \rho(B_{\phi}(\bB\z^*, \bB\z_t) - B_{\phi}(\bB\z^*, \c -\bA\x_{t+1}) - B_{\phi}(\c -\bA\x_{t+1}, \bB\z_t)) \nonumber \\
&+\rho_{\x}(B_{\varphi_{\x}}(\x^*, \x_t) - B_{\varphi_{\x}}(\x^*, \x_{t+1}) - B_{\varphi_{\x}}(\x_{t+1},\x_t))
~.
\end{align}
where the last equality uses the three point property of Bregman divergence, i.e.,
\begin{align}\label{eq:threebregman}
\langle \nabla\phi(\u) - \nabla\phi(\v), \w - \u \rangle = B_{\phi}(\w,\v) - B_{\phi}(\w,\u) - B_{\phi}(\u,\v)~.
\end{align}

Similarly, using the convexity of $g$ and its subgradient given in~\myref{eq:ggd}, for any $\z$,
\begin{align}\label{eq:gconv}
&g(\z_{t+1})  - g(\z) \nonumber \\
&\leq \langle - \bB^T\{ \y_t + \rho ( \nabla\phi  (\bB\z_{t+1}) - \nabla\phi  (\c - \bA\x_{t+1})\} - \rho_{\z}(\nabla\varphi_{\z}(\z_{t+1}) - \nabla\varphi_{\z}(\z_{t})), \z_{t+1} - \z \rangle \nonumber \\
& = -\langle \y_t, \bB(\z_{t+1} - \z) \rangle + \rho  \langle \nabla\phi  (\bB\z_{t+1}) - \nabla\phi(\c - \bA\x_{t+1}) , \bB\z - \bB\z_{t+1}) \rangle \nonumber \\
& - \rho_{\z}\langle \nabla\varphi_{\z}(\z_{t+1}) - \nabla\varphi_{\z}(\z_{t}), \z_{t+1} - \z \rangle \nonumber \\
& = -\langle \y_t, \bB(\z_{t+1} - \z) \rangle + \rho  \left\{ B_{\phi}(\bB\z, \c -\bA\x_{t+1}) - B_{\phi}(\bB\z, \bB\z_{t+1}) -  B_{\phi}(\bB\z_{t+1}, \c -\bA\x_{t+1})\right\} \nonumber \\
& + \rho_{\z}(B_{\varphi_{\z}}(\z, \z_t) - B_{\varphi_{\z}}(\z, \z_{t+1}) - B_{\varphi_{\z}}(\z_{t+1},\z_t))~.
\end{align}
where the last equality uses the three point property of Bregman divergence~\myref{eq:threebregman}.
Set $\z = \z^*$ in~\myref{eq:gconv}. Adding~\myref{eq:fconv} and~\myref{eq:gconv} completes the proof.
\qed

Under Assumption~\ref{asm:badmm}(c), the following lemma shows that~\myref{eq:residual} is bounded by a telescoping series of $D(\w^*,\w_t) - D(\w^*,\w_{t+1})$, where $D(\w^*,\w_t)$ defines the distance from the current iterate $\w_t = (\x_t,\z_t,\y_t)$ to a KKT point $\w^* = (\x^*,\z^*,\y^*)$ as follows:
\begin{align}\label{eq:ht}
D(\w^*,\w_t) = \frac{1}{2\tau\rho}\| \y^* - \y_{t}\|_2^2 + B_{\phi}(\bB\z^*, \bB\z_{t}) + \frac{\rho_{\x}}{\rho} B_{\varphi_{\x}}(\x^*, \x_{t}) + \frac{\rho_{\z}}{\rho} B_{\varphi_{\z}}(\z^*, \z_{t})~.
\end{align}

\begin{lem}\label{lem:RD}
Let the sequence $\{\x_t, \z_t,\y_t\}$ be generated by Bregman ADMM~\myref{eq:badmm_x}-\myref{eq:badmm_y} and $\{\x^*,\z^*,\y^*\}$ satisfying ~\myref{eq:kktx}-\myref{eq:kkty}. Let the Assumption~\ref{asm:badmm} hold. $R(t+1)$ and $D(\w^*,\w_t)$ are defined in~\myref{eq:residual} and~\myref{eq:ht} respectively. Set $\tau \leq (\alpha \sigma  - 2 \gamma)\rho$, where $\sigma = \min\{1, m^{\frac{2}{p}-1}\}$ and $0 < \gamma < \frac{\alpha\sigma}{2}$. Then
\begin{align}\label{eq:case1}
& R(t+1) \leq D(\w^*,\w_t) - D(\w^*,\w_{t+1})~.
\end{align}
\end{lem}
\proof
Assume $\{\x^*,\y^*\}$ satisfies~\myref{eq:kktx}. Since $f$ is convex, then
\begin{align}
f(\x^*) -  f(\x_{t+1}) & \leq - \langle \bA^T\y^*, \x^* - \x_{t+1} \rangle = - \langle \y^*, \bA\x^* - \bA\x_{t+1} \rangle~.
\end{align}
Similarly, for convex function $g$ and $\{\z^*,\y^*\}$ satisfying~\myref{eq:kktz}, we have
\begin{align}
g(\z^*) - g(\z_{t+1}) & \leq - \langle \bB^T\y^*, \z^* - \z_{t+1} \rangle = - \langle \y^*, \bB\z^* - \bB\z_{t+1} \rangle ~.
\end{align}
Adding them together and using the fact that $\bA\x^* + \bB\z^* = \c$, we have
\begin{align}\label{eq:fgopt}
& f(\x^*)+g(\z^*) - (f(\x_{t+1})+g(\z_{t+1})) \leq  \langle \y^*, \bA\x_{t+1} + \bB\z_{t+1} - \c \rangle~.
\end{align}
Adding \myref{eq:fgopt} and \myref{eq:fg} together yields
\begin{align}\label{eq:xyzopt0}
&0 \leq \langle \y^* - \y_t, \bA\x_{t+1} + \bB\z_{t+1} - \c \rangle - \rho ( B_{\phi}(\c -\bA\x_{t+1}, \bB\z_t) +  B_{\phi}(\bB\z_{t+1}, \c -\bA\x_{t+1}))\nonumber \\
&+ \rho (B_{\phi}(\bB\z^*, \bB\z_t) - B_{\phi}(\bB\z^*, \bB\z_{t+1}))+\rho_{\x}(B_{\varphi_{\x}}(\x^*, \x_t) - B_{\varphi_{\x}}(\x^*, \x_{t+1}) - B_{\varphi_{\x}}(\x_{t+1},\x_t))\nonumber \\
& + \rho_{\z}(B_{\varphi_{\z}}(\z^*, \z_t) - B_{\varphi_{\z}}(\z^*, \z_{t+1}) - B_{\varphi_{\z}}(\z_{t+1},\z_t))
~.
\end{align}
Using $\bA\x_{t+1} + \bB\z_{t+1} - \c = \frac{1}{\tau}(\y_{t+1} - \y_t)$, the first term can be rewritten as
\begin{align}
&\langle \y^* - \y_t, \bA\x_{t+1} + \bB\z_{t+1} - \c \rangle = \frac{1}{\tau}\langle \y^* - \y_{t}, \y_{t+1} - \y_t \rangle \nonumber \\
& = \frac{1}{2\tau}\left( \| \y^* - \y_{t}\|_2^2 - \| \y^* - \y_{t+1}\|_2^2 + \| \y_{t+1} - \y_t\|_2^2 \right) \nonumber \\
& = \frac{1}{2\tau}\left( \| \y^* - \y_{t}\|_2^2 - \| \y^* - \y_{t+1}\|_2^2\right) + \frac{\tau}{2}\| \bA\x_{t+1} + \bB\z_{t+1} - \c\|_2^2 ~.
\end{align}
Plugging into~\myref{eq:xyzopt0} and rearranging the terms, we have
\begin{align}
& \frac{1}{2\tau}\left( \| \y^* - \y_{t}\|_2^2 - \| \y^* - \y_{t+1}\|_2^2\right) + \rho (B_{\phi}(\bB\z^*, \bB\z_t) - B_{\phi}(\bB\z^*, \bB\z_{t+1}))  \nonumber \\
&\rho_{\x} (B_{\varphi_{\x}}(\x^*, \x_t) - B_{\varphi_{\x}}(\x^*, \x_{t+1})) + \rho_{\z}(B_{\varphi_{\z}}(\z^*, \z_t) - B_{\varphi_{\z}}(\z^*, \z_{t+1})) \nonumber \\
& \geq \rho_{\x}B_{\varphi_{\x}}(\x_{t+1},\x_t) + \rho_{\z}B_{\varphi_{\z}}(\z_{t+1},\z_t) + \rho B_{\phi}(\c -\bA\x_{t+1}, \bB\z_t) \nonumber \\
& +\rho B_{\phi}(\bB\z_{t+1}, \c -\bA\x_{t+1}) -  \frac{\tau}{2}\| \bA\x_{t+1} + \bB\z_{t+1} - \c\|_2^2~.
\end{align}
Dividing both sides by $\rho$ and letting $R(t+1)$ and $D(\w^*,\w_t)$ be defined in~\myref{eq:residual} and~\myref{eq:ht} respectively, we have
\begin{align}\label{eq:xyzopt}
\!\!\!&D(\w^*,\w_t) \!-\! D(\w^*,\w_{t+1})\! \geq\! R(t+1) \!+\! B_{\phi}(\bB\z_{t+1}, \c -\bA\x_{t+1}) \!-\!  (\frac{\tau}{2\rho}+\gamma)\| \bA\x_{t+1} + \bB\z_{t+1} - \c\|_2^2 \nonumber \\
&\quad \quad\quad\quad\geq R(t+1) +\frac{\alpha}{2} \|\bA\x_{t+1} + \bB\z_{t+1} - \c\|_p^2 -  (\frac{\tau}{2\rho}+\gamma)\| \bA\x_{t+1} + \bB\z_{t+1} - \c\|_2^2
~,
\end{align}
where the last inequality uses the Assumption~\ref{asm:badmm}(c). 

If $0 < p \leq 2$, $\|\u\|_p \geq \|\u\|_2$. Set $\frac{\alpha}{2} \geq \frac{\tau}{2\rho}+\gamma$ in~\myref{eq:xyzopt}, i.e., $\tau \leq (\alpha - 2\gamma) \rho$. We can always find a $\gamma < \frac{\alpha}{2}$, thus ~\myref{eq:case1} follows.

If $p > 2$, $\|\u\|_2 \leq m^{\frac{1}{2} - \frac{1}{p}}\|\u\|_p$ for any $\u \in \R^{m\times 1}$, so $\|\u\|_p^2 \geq m^{\frac{2}{p}-1}\|\u\|_2^2$. In~\myref{eq:xyzopt}, set $\frac{\alpha}{2}m^{\frac{2}{p}-1} \geq \frac{\tau}{2\rho}+\gamma$, i.e., $\tau \leq (\alpha m^{\frac{2}{p}-1} - 2\gamma) \rho$. As long as $\gamma < \frac{\alpha}{2}m^{\frac{2}{p}-1}$, we have ~\myref{eq:case1}.
\qed

\begin{remark}\label{remark1}
(a) If $0< p \leq 2$, then $\sigma = 1$ and $\tau \leq (\alpha - 2\gamma)\rho$. The case that $0 < p \leq 2$ includes two widely used Bregman divergences, i.e., Euclidean distance and KL divergence. For KL divergence in the unit simplex,  we have $\alpha = 1, p = 1$ in the Assumption~\ref{asm:badmm} (c), i.e., $KL(\u,\v) \geq \frac{1}{2}\|\u - \v\|_1^2$~\cite{Beck03}. 

(b) Since we often set $B_\phi$ to be a quadratic function ($p = 2$), the three special cases in Section 2.1 could choose step size $\tau = (\alpha - 2\gamma)\rho$.

(c) If $p > 2$, the proof requires a sufficiently small step size $\tau$, which may not be needed in practice. 
It would be interesting to see whether we can use a same $\tau = O(\rho)$ for any $p > 0$ using other proof techniques. 
\end{remark}

The following theorem establishes the global convergence for BADMM.

\begin{thm} Let the sequence $\{\x_t, \z_t,\y_t\}$ be generated by Bregman ADMM~\myref{eq:badmm_x}-\myref{eq:badmm_y} and $\{\x^*,\z^*,\y^*\}$ satisfying ~\myref{eq:kktx}-\myref{eq:kkty}. Let the Assumption~\ref{asm:badmm} hold and $\tau,\gamma$ satisfy the conditions in Lemma~\ref{lem:RD}. Then $R(t+1)$ converges to zero and $\{\x_t,\z_t,\y_t\}$ converges to a KKT point $\{\x^*,\z^*,\y^*\}$ of~\myref{eq:lag}.
\end{thm}
\proof Since $R(t+1) \geq 0$, \myref{eq:case1} implies
$ D(\w^*,\w_{t+1}) \leq D(\w^*,\w_t)$.
Therefore, $D(\w^*,\w_t)$ is monotonically nonincreasing and $\w_t$ converges to a KKT point $\w^*$.
Summing~\myref{eq:case1} over $t$ from $0$ to $\infty$ yields
\begin{align}\label{eq:case1sum}
&\sum_{t=0}^{\infty} R(t+1) \leq D(\w^*,\w_0)~.
\end{align}
Since $R(t+1) \geq 0$, $R(t+1)\rightarrow 0$ as $t\rightarrow \infty$, which completes the proof.
\qed

The following theorem establishs a $O(1/T)$ convergence rate for the objective and residual of constraints in an ergodic sense.
\begin{thm} \label{thm:badmm_rate}
Let the sequences $\{\x_t, \z_t,\y_t\}$ be generated by Bregman ADMM~\myref{eq:badmm_x},\myref{eq:badmm_z},\myref{eq:badmm_y} and $\y_0 = \mathbf{0}$. Let $\bar{\x}_T = \frac{1}{T}\sum_{t=1}^T \x_t$, $\bar{\z}_T = \frac{1}{T}\sum_{t=1}^T \z_t$. Set $\tau \leq (\alpha \sigma  - 2 \gamma)\rho$, where $\sigma = \min\{1, m^{\frac{2}{p}-1}\}$ and $0 < \gamma < \frac{\alpha\sigma}{2}$. For any $(\x^*,\z^*,\y^*)$ satisfying KKT conditions~\myref{eq:kktx}-\myref{eq:kkty}, we have
\begin{align}
&f(\bar{\x}_{T}) + g(\bar{\z}_{T}) - (f(\x^*) + g(\z^*))
 \leq \frac{D_1}{T}~,\label{eq:fgrate1} \\
&\|\bA\bar{\x}_{T} + \bB\bar{\z}_{T} - \c \|_2^2 \leq \frac{D(\w^*,\w_0)}{\gamma T}~,\label{eq:cstrate1}
\end{align}
where $D_1 = \rho B_{\phi}(\bB\z^*, \bB\z_0)+ \rho_{\x}B_{\varphi_{\x}}(\x^*, \x_0) + \rho_{\z}B_{\varphi_{\z}}(\z^*, \z_0)$.
\end{thm}
\proof Using~\myref{eq:badmm_y}, we have
\begin{align}
-\langle \y_t, \bA\x_{t+1} + \bB\z_{t+1} - \c \rangle &= -\frac{1}{\tau}\langle \y_t, \y_{t+1} - \y_t \rangle \nonumber \\
&= - \frac{1}{2\tau}(\|\y_{t+1}\|_2^2 - \|\y_t\|_2^2 - \| \y_{t+1} - \y_t\|_2^2) \nonumber \\
& = \frac{1}{2\tau}(\|\y_t\|_2^2 - \|\y_{t+1}\|_2^2) + \frac{\tau}{2} \| \bA\x_{t+1} + \bB\z_{t+1} - \c\|_2^2~.
\end{align}
Plugging into~\myref{eq:fg} and ignoring some negative terms yield
\begin{align}\label{eq:fgappen}
&f(\x_{t+1})+ g(\z_{t+1}) - (f(\x^*)+g(\z^*)) \nonumber \\
& \leq  \frac{1}{2\tau}(\|\y_t\|_2^2 - \|\y_{t+1}\|_2^2) + \rho (B_{\phi}(\bB\z^*, \bB\z_t) - B_{\phi}(\bB\z^*, \bB\z_{t+1}))+\rho_{\x}(B_{\varphi_{\x}}(\x^*, \x_t) - B_{\varphi_{\x}}(\x^*, \x_{t+1}))\nonumber \\
& + \rho_{\z}(B_{\varphi_{\z}}(\z^*, \z_t) - B_{\varphi_{\z}}(\z^*, \z_{t+1})) - \rho B_{\phi}(\bB\z_{t+1}, \c -\bA\x_{t+1})+ \frac{\tau}{2} \| \bA\x_{t+1} + \bB\z_{t+1} - \c\|_2^2~.
\end{align}
Assume $B_{\phi}(\bB\z_{t+1}, \c -\bA\x_{t+1}) \geq \frac{\alpha}{2} \| \bA\x_{t+1} + \bB\z_{t+1} - \c \|_p^2$. If $0 < p \leq 2$, using $\|\u\|_p \leq \|\u\|_2$,
\begin{align*}
- \rho B_{\phi}(\bB\z_{t+1}, \c -\bA\x_{t+1}) + \frac{\tau}{2} \| \bA\x_{t+1} + \bB\z_{t+1} - \c\|_2^2 \leq -\frac{\alpha\rho - \tau}{2}\| \bA\x_{t+1} + \bB\z_{t+1} - \c\|_2^2~.
\end{align*}
Setting $\tau \leq (\alpha-2\gamma)\rho$, the last two terms on the right hand side of ~\myref{eq:fgappen} can be removed.

If $p > 2$, $\|\u\|_2 \leq m^{\frac{1}{2} - \frac{1}{p}}\|\u\|_p$ for any $\u \in \R^{m\times 1}$, so $\|\u\|_p^2 \geq m^{\frac{2}{p}-1}\|\u\|_2^2$. Then
\begin{align*}
- \rho B_{\phi}(\bB\z_{t+1}, \c -\bA\x_{t+1})+ \frac{\tau}{2} \| \bA\x_{t+1} + \bB\z_{t+1} - \c\|_2^2 \leq -\frac{\alpha\rho m^{\frac{2}{p}-1} - \tau}{2}\| \bA\x_{t+1} + \bB\z_{t+1} - \c\|_2^2~.
\end{align*}
Setting $\tau \leq (\alpha m^{\frac{2}{p}-1} - 2\gamma)\rho$, the last two terms on the right hand side of ~\myref{eq:fgappen} can be removed. 
Summing over $t$ from $0$ to $T-1$, we have the following telescoping sum
\begin{align}
&\sum_{t=0}^{T-1}\left[f(\x_{t+1}) + g(\z_{t+1}) - (f(\x^*) + g(\z^*))\right] \nonumber \\
&\leq \frac{1}{2\tau}\|\y_0\|_2^2 + \rho B_{\phi}(\bB\z^*, \bB\z_0)+ \rho_{\x}B_{\varphi_{\x}}(\x^*, \x_0) + \rho_{\z}(B_{\varphi_{\z}}(\z^*, \z_0) \nonumber \\
&= \rho B_{\phi}(\bB\z^*, \bB\z_0)+ \rho_{\x}B_{\varphi_{\x}}(\x^*, \x_0) + \rho_{\z}(B_{\varphi_{\z}}(\z^*, \z_0)~.
\end{align}
Dividing both sides by $T$ and applying the Jensen's inequality gives~\myref{eq:fgrate1}.

Dividing both sides of~\myref{eq:case1sum} by $T$ and applying the Jensen's inequality yield~\myref{eq:cstrate1}.
\qed

We consider one special case of BADMM which could outperform ADMM. Assume $\bB = \bI$ and $\cX, \cZ$ are the unit simplex. Let $B_{\phi}$ be the KL divergence. 
For $\z \in \R^{n_2\times 1}$, we have
\begin{align}
B_{\phi}(\z^*, \z_0) = \sum_{i=1}^{n_2} z_i^* \ln \frac{z_i^*}{z_{i,0}} = \sum_{i=1}^{n_2} z_i^* \ln z_i^* + \ln n_2 \leq \ln n_2~.
\end{align}
Similarly, if $\rho_{\x} > 0$, by choosing $\x_0 = \e/n_2$, $B_{\varphi_\x}(\x^*, \x_0) \leq \ln n_1$.
Setting $\alpha = 1, \sigma = 1$ and $\gamma = \frac{1}{4}$ in Theorem~\ref{thm:badmm_rate} yields the following result:
\begin{cor} \label{cor:badmm_rate_KL}
Let the sequences $\{\x_t, \z_t,\y_t\}$ be generated by Bregman ADMM~\myref{eq:badmm_x},\myref{eq:badmm_z},\myref{eq:badmm_y} and $\y_0 = \mathbf{0}$. Assume $\bB = \bI$, and $\cX$ and $\cZ$ is the unit simplex. Let $B_{\phi}, B_{\varphi_\x}, B_{\varphi_\z}$ be KL divergence.
Let $\bar{\x}_T = \frac{1}{T}\sum_{t=1}^T \x_t$, $\bar{\z}_T = \frac{1}{T}\sum_{t=1}^T \z_t$. Set $\tau = \frac{3\rho}{4}$. For any $(\x^*,\z^*,\y^*)$ satisfying KKT conditions~\myref{eq:kkt}, we have
\begin{align}
& f(\bar{\x}_{T}) + g(\bar{\z}_{T}) - (f(\x^*) + g(\z^*)) \leq \frac{\rho\ln n_2  +  \rho_\x \ln n_1 + \rho_{\z} \ln n_2}{T}~,\label{eq:fgrate1} \\
&\|\bA\bar{\x}_{T} + \bB\bar{\z}_{T} - \c \|_2^2 \leq \frac{ \frac{2}{\tau\rho}\| \y^* \!-\! \y_{0}\|_2^2 + 4\ln n_2 + \frac{4\rho_{\x}}{\rho} \ln n_1 \!+\! \frac{4\rho_{\z}}{\rho} \ln n_2 }{ T}~,\label{eq:cstrate1}
\end{align}
\end{cor}

\begin{remark}
(a)  In~\cite{Beck03}, it shows that MDA yields a smilar $O(\ln n)$ bound where $n$ is dimensionality of the problem. If the diminishing step size of MDA is propotional to $\sqrt{\ln n}$, the bound is $O(\sqrt{\ln n})$. Therefore, MDA can outperform the gradient method by a factor $O((n/\ln n)^{1/2})$.

(b) With constant step size, BADMM outperforms ADMM by a factor $O(n/\ln n)$ in an ergodic sense.
\end{remark}

\section{Experimental Results}

In this section,  we use BADMM to solve the mass transportation problem~\cite{cock41:mt}:
\begin{align}\label{eq:mt}
\min&~~~ \langle \bC, \bX\rangle \quad \text{s.t.}\quad  \bX\e = \a ,  \bX^T\e = \b, \bX \geq 0~.
\end{align}
where $\langle \bC,\bX\rangle$ denotes $\tr(\bC^T\bX)$, $\bC \in \R^{m\times n}$ is a cost matrix, $\e$ is a column vector of ones. 
~\myref{eq:mt} is called the assignment problem and can be solved exactly by the Hungarian method~\cite{kuhn55:hungarian}. The mass transportation problem~\myref{eq:mt} is a linear program and thus can be solved by the simplex method.

We now show that~\myref{eq:mt} can be solved by ADMM and BADMM. We first introduce a variable $\bZ$ to split the constraints into two simplex such that $\bd\Delta_{\x} = \{\bX | \bX \geq 0, \bX\e = \a \}$ and $\bd\Delta_{\z} = \{\bZ | \bZ \geq 0, \bZ^T\e = \b \}$.
\myref{eq:mt} can be rewritten in the following ADMM form:
\begin{align}\label{eq:mt_admm}
\min~~~ \langle \bC, \bX\rangle \quad\text{s.t.}\quad  \bX\in \bd\Delta_{\x}, \bZ \in \bd\Delta_{\z},  \bX = \bZ~.
\end{align}
\myref{eq:mt_admm} can be solved by ADMM which requires the Euclidean projection onto the simplex $\bd\Delta_{\x}$ and $\bd\Delta_{\z}$, although  the projection can be done efficiently~\cite{Duchi08}.
We use BADMM to solve~\myref{eq:mt_admm}:
\begin{align}
&\bX^{t+1} = \underset{\bX \in \bd\Delta_{\x}}{\argmin} \langle \bC, \bX \rangle + \langle \bY^t, \bX \rangle + \rho \text{KL} (\bX, \bZ^t )~,\label{eq:mt_badmm_x} \\
&\bZ^{t+1} =  \underset{\bZ \in \bd\Delta_{\z}}{\argmin} \langle \bY^t, - \bZ \rangle + \rho \text{KL} (\bZ, \bX^{t+1}) ~, \label{eq:mt_badmm_z}\\
&\bY^{t+1} = \bY^t + \rho(\bX^{t+1} - \bZ^{t+1})\label{eq:mt_badmm_y}~.
\end{align}
Both~\myref{eq:mt_badmm_x} and~\myref{eq:mt_badmm_z} have closed-form solutions, i.e.,
\begin{align} \label{eq:mult}
X_{ij}^{t+1} &= \frac{Z_{ij}^{t}\exp(-\frac{C_{ij}+Y^t_{ij} }{\rho})}{\sum_{j=1}^n Z_{ij}^{t}\exp(-\frac{C_{ij}+Y^t_{ij}}{\rho})}a_i~, \quad Z_{ij}^{t+1} = \frac{X_{ij}^{t+1}\exp(\frac{ Y^t_{ij}}{\rho})}{\sum_{i=1}^m X_{ij}^{t+1}\exp(\frac{Y^t_{ij}}{\rho})}b_j
\end{align}
which are exponentiated graident updates and can be done in $O(mn)$. Besides the sum operation which can be done in $O(\log(n))$, ~\myref{eq:mult} amounts to elementwise operation and thus can be done in parallel. According to Corollary~\ref{cor:badmm_rate_KL}, BADMM can be faster than ADMM by a factor of $O(n/\log(n))$.

  \begin{figure*}[thb]
   \vspace*{-4mm}
  \subfigure[ m = n = 1000 ]{\label{fig:err_time}
  \includegraphics[width= 45mm,height = 28mm]{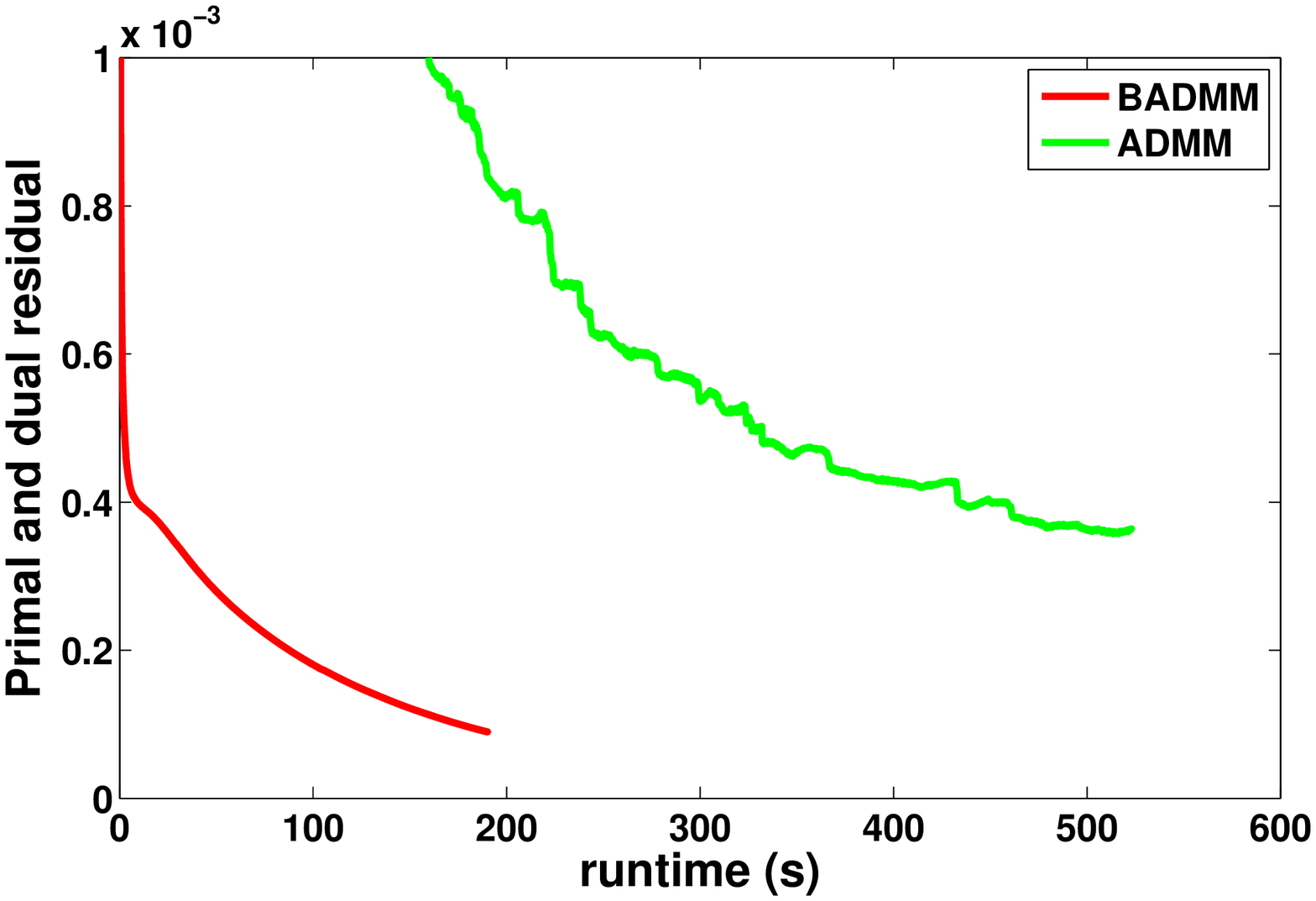}}
  \subfigure[m = n = 2000]{\label{fig:obj_iter}
  \includegraphics[width = 45mm,height = 28mm]{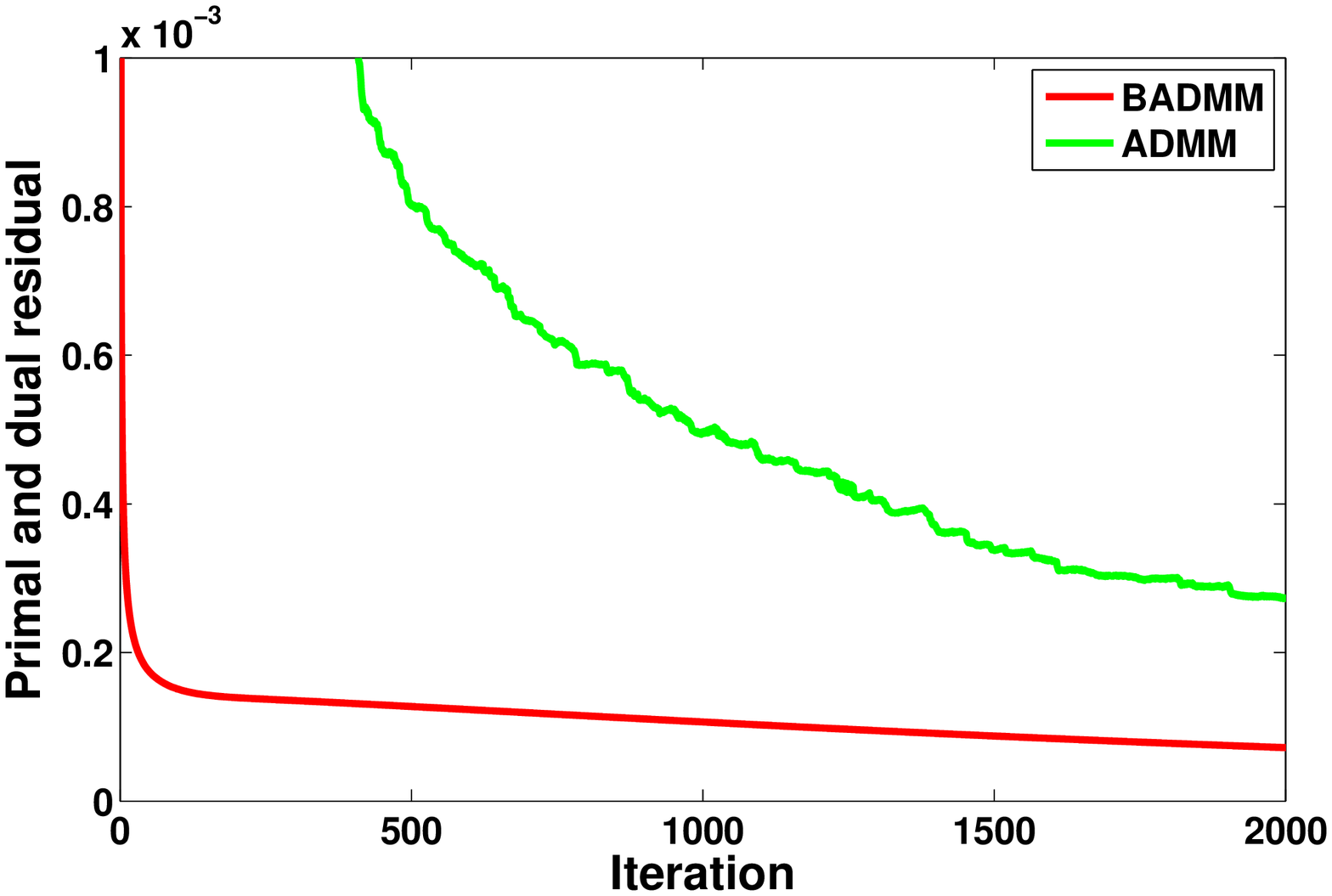}}
  \subfigure[m = n = 4000 ]{\label{fig:obj_time}
  \includegraphics[width = 45mm,height = 28mm]{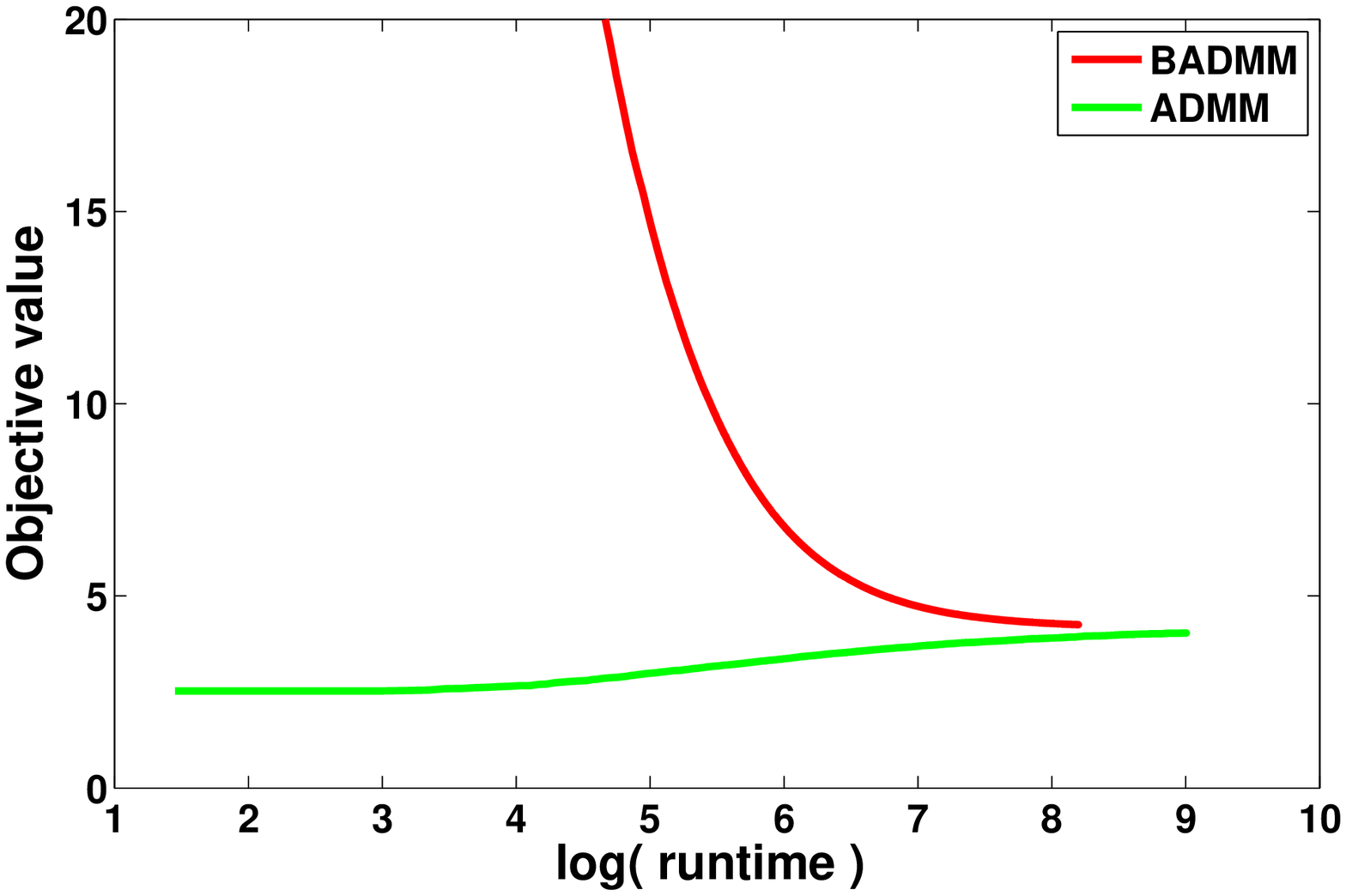}}
  \vspace*{-5mm}
  \caption{ Comparison BADMM and ADMM. BADMM converges faster than ADMM.}\label{fig:badmm_admm} 
   \vspace*{-1mm}
  \end{figure*}

We compare BADMM with ADMM and a highly optimized commercial linear programming solvers on the mass transportation problem~\myref{eq:mt} when $m = n$ and $\a = \b = \e$.  $\bC$ is randomly generated from the uniform distribution. They run 5 times and the average is reported.  We choose the 'best' parameter for BADMM ($\rho = 0.001$) and ADMM ($\rho = 0.001$). The stopping condition is either when the number of iterations exceeds 2000 or when the primal-dual residual is less than $10^{-4}$.

\textbf{BADMM vs ADMM:} Figure~\ref{fig:badmm_admm} compares BADMM and ADMM with different dimensions $n = \{ 1000,2000, 4000 \}$ running on a single CPU. Figure~\ref{fig:err_time} plots the primal and dual residual against the runtime when the dimension is $1000$, and Figure~\ref{fig:obj_iter} plots the convergence of primal and dual residual over iteration when the dimension is $2000$. BADMM converges faster than ADMM.
Figure~\ref{fig:obj_time} plots the convergence of objective value against the log of runtime. BADMM converges faster than ADMM even when the initial point is further from the optimum.


\textbf{BADMM vs Gurobi:}  Gurobi\footnote{\url{http://www.gurobi.com/}} is a highly optimized commercial software where linear programming solvers have been efficiently implemented.
We run Gurobi on two settings: a Mac laptop with 6G memory and a server with 86G memory, respectively. 
For comparison, BADMM is run in parallel on a Tesla M2070 GPU with 5G memory and 448 cores\footnote{GPU code is available on \url{http://www-users.cs.umn.edu/~huwang/badmm_mt.zip}}. 
We experiment with large scale problems and use $m = n = \{1, 5,10,15\}\times 2^{10}$.
Table 1 shows the runtime and the objective values of BADMM and Gurobi, where a `-' indicates the algorithm did not terminate.   
In spite of Gurobi being one of the most optimized LP solvers, BADMM running in parallel is several times faster than Gurobi. In fact, for larger values of $n$, Gurobi did not terminate even on the 86G server, whereas BADMM was efficient even with just 5G memory! The complexity of most LP solvers in Gurobi is $O(n^3)$ and can become slow as $n$ increases, especially at the scales we consider. Moreover,
the memory consumption of Gurobi increases rapidly with the increase of $n$.  When $n = 5\times 2^{10}$, the memory required by Gurobi surpassed the memory in the laptop, leading to the rapid increase of time.  A similar situation was also observed in the server with 86G when $n = 10\times 2^{10}$.  In contrast, the memory required by BADMM is $O(n^2)$---even when $n = 15\times 2^{10}$ (more than 0.2 billion parameters), BADMM can still run on a single GPU with only 5G memory.

The results clearly illustrate the promise of BADMM. With more careful implementation and code optimization, BADMM has the potential to solve large scale problems efficiently in parallel with small memory foot-print.



 \begin{table}
\vspace*{-2mm}
 \caption{ Comparison of BADMM (GPU) with Gurobi}
 \hspace*{-0mm}
 \centering
 \begin{tabular}{|c|c|c|c|c|c|c|}
 \hline
  m=n & \multicolumn{2}{|c|}{Gurobi (Laptop)} & \multicolumn{2}{|c|}{Gurobi (Server)} &  \multicolumn{2}{|c|}{BADMM (GPU)}  \\
  \cline{2-7}
  & time & objective & time & objective & time & objective \\
 \hline
$2^{10}$ &  4.22 & 1.69   & 2.66 & 1.69 & 0.54 & 1.69\\
 \hline
$5\times 2^{10}$  &  377.14 & 1.61  &  92.89 & 1.61 & 22.15 & 1.61 \\
 \hline
$10\times 2^{10}$ & - & - & 1235.34 & 1.65 & 117.75  & 1.65\\
 \hline
 $15\times 2^{10}$ & - & - & - & - &  303.54 & 1.63\\
\hline
 \end{tabular}
 \end{table}\label{tbl:time}

\section{Conclusions}
In this paper, we generalized the alternating direction method of multipliers (ADMM)
to Bregman ADMM, similar to how mirror descent generalizes gradient descent. BADMM defines a unified framework for ADMM, generalized ADMM, inexact ADMM and Bethe ADMM.
The global convergence and the $O(1/T)$ iteration complexity of BADMM are also established. In some cases, BADMM is faster than ADMM by a factor of $O(n/\log(n))$.  BADMM can also be faster than highly optimized commercial software in solving linear program of mass transportation problem.

\section*{Acknowledgment}
H.W. and A.B. acknowledge the support of NSF via IIS-0953274, IIS-1029711, IIS- 0916750, IIS-0812183, NASA grant NNX12AQ39A, and the technical support from the University of Minnesota Supercomputing Institute. H.W. acknowledges the support of DDF (2013-2014) from the University of Minnesota. A.B. acknowledges support from IBM and Yahoo.

\bibliographystyle{plain}
\bibliography{long,admm,onlinelearn,all,map,mt}

\appendix

\section{Convergence of BADMM with Time Varying Step Size} \label{appen3}
Under the assumption that $\y_t$ is bounded, the following theorem requires a large step size to establish the convergence of BADMM.

\begin{thm}\label{thm:cvg2}
Let the sequences $\{\x_t, \z_t,\y_t\}$ be generated by Bregman ADMM~\myref{eq:badmm_x}-\myref{eq:badmm_y} and $\{\x^*,\z^*,\y^*\}$ satisfying ~\myref{eq:kktx}-\myref{eq:kkty}. Let the Assumption~\ref{asm:badmm} hold and $\|\y_t\|_2 \leq D_\y$. Setting $\rho_\x = \rho_\z = c_1\sqrt{T}, \tau = c_2\sqrt{T}$ and $\rho = \sqrt{T}$ for some positive constant $c_1,c_2$, then $R(t+1)$ converges to zero.
\end{thm}
\proof
Assuming $\|\y_t\|_2 \leq D_\y$ and using~\myref{eq:badmm_y}, we have
\begin{align}
\|\bA\x_{t+1} + \bB\z_{t+1} - \c\|_2^2 = \frac{1}{\tau^2} \| \y_{t+1} - \y_t \|_2^2 \leq \frac{2}{\tau^2} (\| \y_{t+1} \|_2^2+\| \y_{t} \|_2^2)\leq \frac{4D_\y^2}{\tau^2}~.
\end{align}
Plugging into~\myref{eq:xyzopt} and rearranging the terms yields
\begin{align}\label{eq:case2}
& R(t+1) \leq D(\w^*,\w_t) - D(\w^*,\w_{t+1}) + (\frac{\tau}{2\rho}+\gamma)\frac{4D_\y^2}{\tau^2}~.
\end{align}
Setting $\rho_\x = \rho_\z = c_1\sqrt{T}, \tau = c_2\sqrt{T}$ and $\rho = \sqrt{T}$ for some positive constant $c_1,c_2$, we have
\begin{align}
R(t+1) = c_1 B_{\varphi_{\x}}(\x_{t+1},\x_t) + c_1 B_{\varphi_{\z}}(\z_{t+1},\z_t) + B_{\phi}(\c -\bA\x_{t+1}, \bB\z_t) +\gamma\| \bA\x_{t+1} + \bB\z_{t+1} - \c\|_2^2~,
\end{align}
Summing~\myref{eq:case2} over $t$ from $0$ to $T-1$, we have the following telescoping sum
\begin{align}\label{eq:case2sum}
\sum_{t=0}^{T-1} R(t+1) &\leq D(\w^*,\w_0) + \sum_{t=0}^{T-1}(\frac{\tau}{2\rho}+\gamma)\frac{4D_\y^2}{\tau^2} = D(\w^*,\w_0) + \frac{4(c_2/2+\gamma)D_\y^2}{c_2^2}~.
\end{align}
Therefore, $R(t+1) \rightarrow 0$ as $t\rightarrow \infty$.
\qed

The following theorem establishs the convergence rate for the objective and residual of constraints in an ergodic sense.
\begin{thm} \label{thm:badmm_rate2}
Let the sequences $\{\x_t, \z_t,\y_t\}$ be generated by Bregman ADMM~\myref{eq:badmm_x}-\myref{eq:badmm_y}. Let $\bar{\x}_T = \frac{1}{T}\sum_{t=1}^T \x_t, \bar{\z}_T = \frac{1}{T}\sum_{t=1}^T \z_t$.
Let the Assumption~\ref{asm:badmm} hold and $\|\y_t\|_2 \leq D_\y$. Set $\rho_\x = \rho_\z = c_1\sqrt{T}, \tau = c_2\sqrt{T}, \rho = \sqrt{T}$ for some positive constants $c_1,c_2$. For any $(\x^*,\z^*,\y^*)$ satisfying KKT conditions~\myref{eq:kktx}-\myref{eq:kkty}, we have
\begin{align}
&f(\bar{\x}_{T}) + g(\bar{\z}_{T}) - (f(\x^*) + g(\z^*))
 \leq \frac{2D_\y^2}{c_2\sqrt{T}} + \frac{\|\y_0\|_2^2}{2c_2T\sqrt{T}} + \frac{D_2}{\sqrt{T}}~,\label{eq:fgrate2} \\
&\|\bA\bar{\x}_{T} + \bB\bar{\z}_{T} - \c \|_2^2 \leq \frac{D(\w^*,\w_0)}{\gamma T} + \frac{4(c_2/2+\gamma)D_\y^2}{\gamma c_2^2T}~,\label{eq:cstrate2}
\end{align}
where $D_2 = B_{\phi}(\bB\z^*, \bB\z_0) + c_1(B_{\varphi_{\x}}(\x^*, \x_0) + B_{\varphi_{\z}}(\z^*, \z_0))$.
\end{thm}
\proof
Assuming $\|\y_t\|_2 \leq D^2_\y$ and using~\myref{eq:badmm_y}, we have
\begin{align}
-\langle \y_t, \bA\x_{t+1} + \bB\z_{t+1} - \c \rangle = -\frac{1}{\tau} \langle \y_t, \y_{t+1} - \y_t \rangle \leq \frac{1}{\tau}(\|\y_t\|_2^2 + \|\y_t\|_2*\|\y_{t+1}\|_2) \leq \frac{2D^2_\y}{\tau}~.
\end{align}
Plugging into~\myref{eq:fg} and ignoring some negative terms yield
\begin{align}
&f(\x_{t+1})+ g(\z_{t+1}) - (f(\x^*)+g(\z^*)) \nonumber \\
& \leq  \frac{2D^2_\y}{\tau} + \rho (B_{\phi}(\bB\z^*, \bB\z_t) - B_{\phi}(\bB\z^*, \bB\z_{t+1}))+\rho_{\x}(B_{\varphi_{\x}}(\x^*, \x_t) - B_{\varphi_{\x}}(\x^*, \x_{t+1}))\nonumber \\
& + \rho_{\z}(B_{\varphi_{\z}}(\z^*, \z_t) - B_{\varphi_{\z}}(\z^*, \z_{t+1}))~.
\end{align}
Summing over $t$ from $0$ to $T-1$, we have the following telescoping sum
\begin{align*}
&\sum_{t=0}^{T-1}\left[f(\x_{t+1}) + g(\z_{t+1}) - (f(\x^*) + g(\z^*))\right] \nonumber \\
&\leq \sum_{t=0}^{T-1}\frac{2D^2_\y}{\tau} + \frac{1}{2\tau}\|\y_0\|_2^2 + \rho B_{\phi}(\bB\z^*, \bB\z_0)+ \rho_{\x}B_{\varphi_{\x}}(\x^*, \x_0) + \rho_{\z}B_{\varphi_{\z}}(\z^*, \z_0)~.
\end{align*}
Setting $\rho_\x = \rho_\z = c_1\sqrt{T}, \tau = c_2\sqrt{T}, \rho = \sqrt{T}$, dividing both sides by $T$ and applying the Jensen's inequality yield ~\myref{eq:fgrate2}.

Dividing both sides of~\myref{eq:case2sum} by $T$ and applying the Jesen's inequality yield~\myref{eq:cstrate2}.
\qed

\end{document}